\documentclass[a4paper,11pt]{article}
\usepackage[T1]{fontenc} 
\usepackage[cp1250]{inputenc}
\usepackage{lmodern}

\usepackage[english]{babel}
\usepackage{latexsym}
\usepackage{amsmath}
\usepackage{amsthm}
\usepackage{amssymb}
\usepackage{indentfirst}
\usepackage{amsfonts}
\usepackage{stackrel}
\usepackage{dsfont}
\usepackage{tikz}
\usepackage{slashbox}
\usepackage[mathscr]{eucal}
\usepackage{authblk}

\newtheorem{thm}{Theorem}
\newtheorem{lemma}[thm]{Lemma}

\newtheorem{defin}{Definition}
 
\newenvironment{pro}{\begin{flushleft} \textbf{Proof}\\* \end{flushleft}}{\hfill\(\blacksquare\) \\ }

\newcommand{\reals}{\mathbb{R}}
\newcommand{\naturals}{\mathbb{N}}

\newcommand{\complex}{\mathbb{C}}

\newcommand{\eps}{\varepsilon}

\newcommand{\conv}{\text{conv}}
\newcommand{\cone}{\mathcal{C}}

\newcommand{\bounded}{\mathcal{M}}

\newcommand{\compact}{\mathcal{N}}

\newcommand{\Ffamily}{\mathcal{F}}
\newcommand{\Gfamily}{\mathcal{G}}

\newcommand{\Hammer}{\mathscr{H}}

\newcommand{\BL}{\mathfrak{B}}

\newcommand{\Addresses}{{
  \bigskip
  \footnotesize

  Krukowski M. (corresponding author), \textsc{Technical Univeristy of \L\'od\'z, \ Institute of Mathematics, \ W\'ol\-cza\'n\-ska 215, \
90-924 \ \L\'od\'z, \ Poland}\par\nopagebreak
  \textit{E-mail address} : \texttt{krukowski.mateusz13@gmail.com}

  \medskip
}}

\begin{document}
\title{Darbo-type theorem for quasimeasure of noncompactness}
\author{Mateusz Krukowski}
\affil{Technical University of \L\'od\'z, Institute of Mathematics, \\ W\'ol\-cza\'n\-ska 215, \
90-924 \ \L\'od\'z, \ Poland}
\maketitle

\begin{abstract}
The paper introduces the concept of quasimeasure of noncompactness. Motivated by the Arzel\`a-Ascoli theorem for $C^b(X,E)$, where $X$ is an Euclidean space and $E$ an arbitrary Banach space, we construct a quasimeasure for this space and study its properties. An analogon for Darbo fixed point theorem is obtained with the additional aid of measure of nonconvexity. The paper ends with possible application in integral equations. We prove that a Hammerstein operator with Carath\'{e}odory kernel and nonlinearity of a certain type has a fixed point. 
\end{abstract}


\section{Introduction}

In 1955, Gabriele Darbo published the paper 'Punti uniti in trasformazioni a codominio non compatto' (comp. \cite{Darbo}), where he introduced what has become known as Darbo theorem. The charming beauty of the theorem lies in the unification of two approaches to fixed point theory. The first one began in 1922 with the paper 'Sur les op\'{e}rations dans les ensembles abstraits et leur application aux \'{e}quations int\'{e}grales' (comp. \cite{Banach}) and today is recognized as Banach contraction principle. Nearly 10 years after Banach's breakthrough (in 1930), Leray Schauder came up with an innovative (topological) way of producing a fixed point (comp. \cite{Schauder}). Today, Schauder fixed point theorem is an indispensable tool in the field of differential equations.  

Year 1930 was the advent of yet another branch of mathematical analysis. Kazimierz Kuratowski in paper 'Sur les espaces complets' introduced the concept of Kuratowski measure of noncompactness. For a metric space $X$, it is a functional defined on bounded sets $B$ via the formula

$$\alpha(B) = \inf \bigg\{ \delta > 0 \ : \ B \subset B_1 \cup \ldots B_n, \ \text{diam}(B_k)\leq \delta, \ 1\leq k\leq n\bigg\}$$

\noindent
With the aid of functional $\alpha$, Kuratowski generalized a well-known Cantor intersection theorem. Over 25 years later (in 1957), a similar measure of noncompactness was introduced and named in the honour of Felix Hausdorff:

$$\beta(B) = \inf \bigg\{ \delta > 0 \ : \ B \subset B(x_1,\delta) \cup \ldots B(x_n,\delta), \  x_1,\ldots,x_n \in X \bigg\}$$

Kuratowski and Hausdorff measure of noncompactnes are equivalent in the sense 

$$\beta(B) \leq \alpha(B) \leq 2\beta(B)$$

\noindent
Moreover, they both resemble the celebrated Hausdorff $\eps$-net theorem, discussed in \cite{Sutherland} on page 146. This is no coincidence, as $\alpha$ is complete in the sense 

\begin{gather}
\alpha(B) = 0 \ \Longleftrightarrow \ B \text{ is relatively compact}
\label{completemeasure}
\end{gather}

\noindent
The same is true for $\beta$. This crucial property of measures of noncompactness made the author realize, that whenever a characterization of compact sets of a space is given, we can try to 'cook up' a suitable object resembling a measure of noncompactness. Following this train of thought, the author has decided to build a quasimeasure of noncompactness based on the Arzel\`a-Ascoli theorem in \cite{KrukowskiPrzeradzki}, which he proved together with Bogdan Przeradzki. We recall this result as theorem \ref{AAforXE}.

The prefix 'quasi-' in quasimeasure of noncompactness is (at this stage) necessary, as not all natural properties of Kuratowski or Hausdorff measure of noncompactness carry over to the constructed object. It may be the case that the author lacks the ability to prove these properties and is forced to work without them. Luckily, the analogon of (\ref{completemeasure}) turns out to be true.

As far as the construction of the paper is concerned, the next subsection establishes the notation which we hold on to throughout the paper. A formal definition of quasimeasure of noncompactness is given and moreover, the concept of measure of nonconvexity is recalled. Section \ref{sectionquasimeasureofnoncompactness} opens with Arzel\`a-Ascoli theorem, proved in \cite{KrukowskiPrzeradzki}. Next, three functionals are considered: $\eta$ which measures pointwise relative compactness, $\omega_0$ which measures equicontinuity and $\chi_0$ which measures $C^b(X,E)$-extension property. The rest of the section is devouted to verifying properties of these three functionals. The culiminating point is theorem \ref{quasioncbxe}, which constructs the desired quasimeasure of noncompactness. 

The section \ref{applications} contains a proof of Darbo-like theorem for quasimeasure of noncompactness. In order to illustrate this theorem at work, we consider a Hammerstein operator with Carath\'{e}odory kernel and look for its fixed point. The paper ends with a digression on what situations are generalized by theorem \ref{fixedpointdarbo}. The evident similarities with Darbo theorem are briefly described. At last, the author admits that the purpose of obtained results lies in unifying mathematical concepts rather than in applications to ceratin kinds of differential or integral equations.

\subsection{Notation and basic definitions}

Throughout the paper, $Y$ denotes a Banach space. $B(c,r)$ and $\overline{B}(c,r)$ are understood to be open and closed balls, respectively, centered at $c$ and of radius $r>0$. The space of bounded and linear operators on Banach space $Y$ is denoted by $\BL(Y)$, while $\Hammer$ stands for Hammerstein operator. Furthermore, we distinguish three families of space $Y$:

\noindent
\begin{description}
	\item[\hspace{0.4cm} $\bounded(Y)$] - the family of bounded subsets of $Y$
	\item[\hspace{0.4cm} $\compact(Y)$] - the family of compact subsets of $Y$
\end{description}

The following definition introduces the concept of quasimeasure of noncompactness. The list of axioms \textbf{(QMN1)} - \textbf{(QMN5)} differs from the classical collection found in \cite{BanasGoebel} (page 11) or \cite{BanasMursaleen} (page 170).

\begin{defin}(quasimeasure of noncompactness)\\
A function $\Omega : \bounded(Y) \rightarrow \reals_+$ is called a quasimeasure of noncompactness in $Y$ if

\begin{description}
	\item[\hspace{0.4cm} (QMN1)] for every $A \in \bounded(Y), \ \Omega(A) = 0$ if and only if $A$ is relatively compact	
	\item[\hspace{0.4cm} (QMN2)] for every $A,B \in \bounded(Y)$ such that $A \subset B$ we have $\Omega(A) \leq \Omega(B)$
	\item[\hspace{0.4cm} (QMN3)] for every $A \in \bounded(Y)$ we have $\Omega(\overline{A}) = \Omega(A)$
	\item[\hspace{0.4cm} (QMN4)] for every $A \in \bounded(Y)$ and $\lambda \in \complex$ we have $\Omega(\lambda A) = |\lambda|\Omega(A)$
	\item[\hspace{0.4cm} (QMN5)] for every finite $A$ and $B \in \bounded(Y)$ such that $A \cap B = \emptyset$ we have $\Omega(A\cup B) = \Omega(B)$
\end{description}

\label{quasimeasureofnoncompactness}
\end{defin}

We also recall the concept of measure of nonconvexity as presented in \cite{Eisenfeld}. It estimates how far a set $A$ lies from its convex hull. The distance $d_H$ is understood to be the Hausdorff distance (\cite{Aliprantis}, page 110, definition 110).

\begin{defin}
Let $(Y,\|\cdot\|_Y)$ be a Banach space. A function $\kappa : \bounded(Y) \rightarrow \reals_+$ given by 

$$\forall_{A \in \bounded(Y)} \ \kappa(A) := d_H(A,\conv(A))$$

\noindent
is called a measure of nonconvexity, where $d_H$ is the Hausdorff distance. 
\label{mofnonconvexity}
\end{defin}

The second part of the paper focuses on integral equations. From this point forward, $X$ denotes an Euclidean space with Lebesgue measure $\lambda$. Moreover, $(E,\|\cdot\|_E)$ denotes a Banach space of arbitrary dimension.

\section{Quasimeasure of noncompactness on $C^b(X,E)$}
\label{sectionquasimeasureofnoncompactness}
Let us recall a theorem (we adjust the notation to this paper), which comes from \cite{KrukowskiPrzeradzki} and characterizes relatively compact families in $C^b(X,E)$.

\begin{thm}
The family $\Ffamily \subset C^b(X,E)$ is relatively compact iff

\begin{description}
	\item[\hspace{0.4cm} (AA1)] $\Ffamily$ is pointwise relatively compact and equicontinuous at every point of $X$
	\item[\hspace{0.4cm} (AA2)] $\Ffamily$ satisfies $C^b(X,E)$-extension property, i.e.
	
	$$\forall_{\eps > 0} \ \exists_{\substack{D \Subset X \\ \delta > 0}} \ \forall_{f,g \in \Ffamily} \ d_{C(D,E)}(f,g) \leq \delta \ \Longrightarrow \ d_{C^b(X,E)}(f,g) \leq \eps $$
\end{description}
\label{AAforXE}
\end{thm}

In light of the above theorem, constructing a quasimeasure of noncompactness, we need to be able to somehow measure three quantities: pointwise relative compactness, pointwise equicontinuity and finally $C^b(X,E)$-extension property. The following definition takes care of pointwise relative compactness.

\begin{defin}
For $\Ffamily \in \bounded(C^b(X,E))$ define

\begin{gather}
\eta(\Ffamily) := \sup_{x \in X} \ \beta(\{ f(x) \ : \ f \in \Ffamily \})
\label{pointwisecompactmeasure}
\end{gather}

\noindent
where $\beta$ is the Hausdorff measure of noncompactness on $E$. 
\end{defin}

In case $X = \reals$, such a functional has been employed in the study of differential equations, for example in \cite{PrzeradzkiAPM}. The functional $\eta$ measures how much a set deviates from being pointwise relatively compact.  The next lemma puts this claim in formal mathematical setting.

\begin{lemma}
The function $\eta : \bounded(C^b(X,E)) \rightarrow \reals_+$ satisfies \emph{\textbf{(QMN2)}-\textbf{(QMN5)}}. Moreover, for $\Ffamily \in \bounded(C^b(X,E))$ we have $\eta(\Ffamily) = 0$ iff $\Ffamily$ is pointwise relatively compact.
\label{pointwisecompactmeasurelemma}
\end{lemma}
\begin{pro}

\textbf{(QMN2)}, \textbf{(QMN4)} and pointwise relative compactness of $A$ iff $\eta(A)=0$ are obvious due to properties of $\beta$ and definition of $\eta$. 

In order to prove \textbf{(QMN3)} observe that if $f_n \in \Ffamily, \ f \in \overline{\Ffamily}$ and $d_{C^b(X,E)}(f_n,f) \rightarrow 0$ then $f_n(x) \rightarrow f(x)$ for every $x \in X$. Hence the inclusion

$$ \overline{\{ f(x) \ : \ f \in \Ffamily \}} \supset \{ f(x) \ : \ f \in \overline{\Ffamily} \} $$
	
\noindent
holds and consequently, by \textbf{(QMN2)} we obtain $\sup_{x \in X} \ \beta(\overline{\{ f(x) \ : \ f \in A \}}) \geq \eta(\overline{A})$. Since $\beta$ is invariant under taking the closure of a set, we are done.

Finally, we need to verify \textbf{(QMN5)}. Observe that if $\Ffamily \in \compact(BC(X,E))$ then by theorem \ref{AAforXE} we have that $\Ffamily$ is pointwise relatively compact. Consequently, $\eta(\Ffamily) = 0$, i.e. 

$$\forall_{x \in X} \ \beta(\{f(x) \ : \ f \in \Ffamily\}) = 0$$

\noindent
which means precisely that $\Ffamily$ is pointwise relatively compact.
\end{pro}

Apart from measuring pointwise relative compactness, we would like to measure the violation of equicontinuity and $C^b(X,E)$-extension property. Further two lemmas introduce sufficient tools for this task.

\begin{lemma}
For every $\Ffamily \in \bounded(C^b(X,E))$ and $x \in X, \ \delta > 0$ we define

\begin{gather}
\begin{split}
	\omega^x(\Ffamily,\delta) &:= \sup_{f \in \Ffamily} \ \sup_{y \in B(x,\delta)} \ \|f(y) - f(x)\| \\
	&\omega_0^x(\Ffamily) := \lim_{\delta \rightarrow 0} \ \omega^x(\Ffamily,\delta) \\
	&\omega_0(\Ffamily) := \sup_{x \in X} \ \omega_0^x(\Ffamily)
\end{split}
\label{continuitymeasure}
\end{gather}

\noindent
The function $\omega_0$ is well-defined and for every $\Ffamily \in \bounded(C^b(X,E))$. Moreover, $\omega_0$ satisfies \emph{\textbf{(QMN2)}-\textbf{(QMN5)}} and $\omega_0(\Ffamily) = 0$ iff $\Ffamily$ is pointwise equicontinuous.
\label{equicontinuitymeasure}
\end{lemma}
\begin{pro}

To prove that $\omega_0$ is well-defined (i.e. the limit $\omega_0^x(\Ffamily)$ exists), it suffices to observe that $\delta \mapsto \omega^x(\Ffamily,\delta)$ is nondecreasing due to the inclusion $B(x,\delta) \subset B(x,\delta')$ for $\delta \leq \delta'$. Furthermore, observe that

\begin{gather*}
\omega_0(\Ffamily) = 0 \ \Longleftrightarrow \ \forall_{x \in X} \ \forall_{\eps > 0} \ \exists_{\delta > 0} \ \omega^x(\Ffamily,\delta) \leq \eps
\end{gather*}

\noindent
which is equivalent to $\Ffamily$ being pointwise equicontinuous at every $x \in X$. 

Properties \textbf{(QMN2)} and \textbf{(QMN4)} are obvious. In order to verify $\textbf{(QMN3)}$ we fix $\eps > 0$ and observe that there exists $f_{\ast} \in \overline{\Ffamily}$ such that 
	
$$\forall_{x \in X} \ \omega^x(\overline{\Ffamily},\delta) \leq \sup_{y \in B(x,\delta)} \ \|f_{\ast}(y) - f_{\ast}(x)\| + \frac{\eps}{3}$$
	
\noindent
Moreover, there exists $g_{\ast} \in \Ffamily$ such that $d_{C^b(X,E)}(f_{\ast},g_{\ast}) < \frac{\eps}{3}$. Then we have
	
$$\omega^x(\overline{\Ffamily},\delta) \leq \sup_{y \in B(x,\delta)} \ \|g_{\ast}(y) - g_{\ast}(x)\| + \eps \leq \omega^x(\Ffamily,\delta) + \eps$$
	
\noindent
We may conclude that $\omega(\overline{\Ffamily}) \leq \omega(\Ffamily)$, which we aimed for. 

Finally, to prove \textbf{(QMN5)} we again use theorem \ref{AAforXE} which implies that if $\Ffamily \in \compact(C^b(X,E))$ then $\omega_0(\Ffamily)=0$. Now it is easy to see that for $\Gfamily \in \bounded(C^b(X,E))$ we have

$$\omega_0(\Gfamily) \leq \omega_0(\Ffamily \cup \Gfamily) \leq \omega_0(\Ffamily) + \omega_0(\Gfamily) = \omega_0(\Gfamily)$$

\noindent
In fact, we proved more than \textbf{(QMN5)}: $\Ffamily$ can be any compact family in $C^b(X,E)$, not merely a finite one.
\end{pro}

\begin{thm}
For every $\Ffamily \in \bounded(C^b(X,E)), \ n\in\naturals$ and $\eps > 0$ we define

\begin{gather}
\chi^n(\Ffamily,\eps) := \sup \ \bigg\{ d_{C^b(X,E)}(f,g) \ : \ f,g \in \Ffamily, \ d_{C(S_n,E)}(f,g) \leq \eps \bigg\}
\label{definchi}
\end{gather}
	
\noindent
where $(S_n)_{n\in\naturals}$ is a saturating sequence for $X$, that is an ascending sequence of compact sets whose union is the whole space $X$. Moreover, 
	
\begin{equation}
	\begin{split}
		&\chi^n_0(\Ffamily) := \lim_{\eps \rightarrow 0} \ \chi^n(\Ffamily,\eps) \\
		&\chi_0(\Ffamily) := \lim_{n \rightarrow \infty} \ \chi^n_0(\Ffamily)
	\end{split}
	\label{limitschi}
\end{equation}

\noindent
The function $\chi_0 : \bounded(C^b(X,E)) \rightarrow \reals_+$ is well-defined, i.e. limits \emph{(\ref{limitschi})} exist and the definition does not depend on the choice of saturating sequence. Moreover, $\chi_0$ satisfies \emph{\textbf{(QMN2)}-\textbf{(QMN5)}} and for every $\Ffamily \in \bounded(C^b(X,E))$ we have $\chi_0(\Ffamily) = 0$ iff $\Ffamily$ satisfies $C^b(X,E)$-extension property.
\label{chilemma}
\end{thm}
\begin{pro}

Throughout the proof, we assume that $\Ffamily \in \bounded(C^b(X,E))$. To show that $\chi_0$ is well defined, observe that $\eps \mapsto \chi^n(A,\eps)$ is nondecreasing and $n \mapsto \chi^n_0(A)$ is nonincreasing due to the implications 

$$d_{BC(S_n,E)}(f,g) \leq \eps \ \Longrightarrow \ d_{BC(S_n,E)}(f,g) \leq \eps'$$

\noindent
for $\eps \leq \eps'$ and 

$$d_{C(S_n,E)}(f,g) \leq d_{C(S_{n'},E)}(f,g)$$

\noindent
for $n \leq n'$.

Moreover, if $(S_n)_{n\in\naturals}$ and $(\tilde{S}_n)_{n\in\naturals}$ are two saturating sequences, then for a fixed $n\in\naturals$, we can find $m_n \geq n$ such that $\tilde{S}_{m_n} \supset S_n$. Hence $\tilde{\chi}^{m_n}(\Ffamily,\eps) \leq \chi^{n}(\Ffamily,\eps)$ and passing to the limits (\ref{limitschi}) we obtain $\tilde{\chi_0}(\Ffamily) \leq \chi_0(\Ffamily)$. Applying the same reasoning in the other direction, we conclude that $\chi_0$ is independent of the choice of saturating sequence.

Using the fact that $n \mapsto \chi^n(\Ffamily,\delta)$ is nonincreasing and $\delta \mapsto \chi^n(\Ffamily,\delta)$ is nondecreasing, the following calculation

\begin{equation}
\chi_0(\Ffamily) = 0 \ \Longleftrightarrow \ \forall_{\eps > 0} \ \exists_{n \in \naturals} \ \forall_{m \geq n} \ \lim_{\delta \rightarrow 0} \chi^m(\Ffamily,\delta) \leq \eps \ \Longleftrightarrow \ \forall_{\eps > 0} \ \exists_{n \in \naturals} \ \exists_{\delta > 0} \ \chi^n(\Ffamily,\delta) \leq \eps   
\label{chiequalszero}
\end{equation}

\noindent
implies that $\chi_0(\Ffamily) = 0$ is equivalent to $\Ffamily$ satisfying $C^b(X,E)$-extension property.

Observe that \textbf{(QMN2)} is obvious due to the properties of supremum. In order to show \textbf{(QMN3)}, we will use $\eps_1, \eps_2 > 0$ as variables (tending to $0$). For (temporarily fixed) $\eps_2$ we choose 
	
\begin{itemize}
	\item $n \in \naturals$ such that $\chi_0^n(\Ffamily) \leq \chi_0(\Ffamily) + \frac{\eps_2}{3}$
	\item $\eps_1 > 0$ such that $\chi^n(\Ffamily,3\eps_1) \leq \chi_0^n(\Ffamily) + \frac{\eps_2}{3}$ and $9 \eps_1 \leq \eps_2$
\end{itemize}
	
\noindent
By (\ref{definchi}), there exist $\overline{f},\overline{g} \in \overline{\Ffamily}$ with $d_{C(S_n,E)}(\overline{f},\overline{g}) \leq \eps_1$ such that
	
\begin{gather}
\chi^n(\overline{\Ffamily},\eps_1) \leq d_{C^b(X,E)}(\overline{f},\overline{g}) + \eps_1
\label{provingclosureproperty}
\end{gather}
	
\noindent
By definition of closure, there exist $f,g \in \Ffamily$ such that $d_{C^b(X,E)}(\overline{f},f)\leq \eps_1$ and $d_{C^b(X,E)}(\overline{g},g) \leq \eps_1$. Hence, the inequality $d_{C(S_n,E)}(\overline{f},\overline{g}) \leq \eps_1$ implies that

\begin{gather}
d_{C(S_n,E)}(f,g) \leq 3\eps_1
\label{trivialestimate}
\end{gather}
				
At last, we have
	
\begin{gather*}
\chi^n(\overline{\Ffamily},\eps_1) \stackrel{(\ref{provingclosureproperty})}{\leq}  d_{C^b(X,E)}(\overline{f},\overline{g}) + \eps_1 \leq d_{C^b(X,E)}(f,g) + 3\eps_1 \stackrel{(\ref{trivialestimate})}{\leq} \chi^n(\Ffamily,3\eps_1) + \eps_2 \\
\stackrel{\text{choice of }\eps_1}{\leq} \chi_0^n(\Ffamily) + \frac{2\eps_2}{3} \stackrel{\text{choice of }n}{\leq} \chi_0(\Ffamily) + \eps_2
\end{gather*}
	
\noindent
Passing to the limits (\ref{limitschi}) and using the fact that $\eps_2$ was chosen arbitrarily, we prove \textbf{(QMN3)}.
	
Next, we verify \textbf{(QMN4)}. We start by observing thta if $\lambda = 0$ there's nothing to prove. If $\lambda \in \complex\backslash\{0\}$, then for $\Ffamily \in \bounded(C^b(X,E))$ we have 

\begin{equation}
\begin{split}
&\chi^n(\lambda \Ffamily,\eps) = \sup\bigg\{ d_{C^b(X,E)}(\lambda f, \lambda g) \ : \ f,g \in \Ffamily, \ d_{C(S_n,E)}(\lambda f, \lambda g) \leq \eps \bigg\} \\
&= |\lambda| \sup\bigg\{ d_{C^b(X,E)}(f,g) \ : \ f,g \in \Ffamily, \ d_{C(S_n,E)}(f,g) \leq \frac{\eps}{|\lambda|} \bigg\} = |\lambda| \chi^n\bigg(\Ffamily,\frac{\eps}{|\lambda|}\bigg) 
\end{split}	
\label{estimateswithlambda}
\end{equation}
		
\noindent
Taking the limits (\ref{limitschi}), we obtain the desired property. 

Last but not least, we show \textbf{(QMN5)}. From this point till the end of the proof, we consider a finite family $\Ffamily$ and $\Gfamily \in \bounded(C^b(X,E))$. Without loss of generality, we may assume that $\Ffamily \cap \overline{\Gfamily} = \emptyset$ due to \textbf{(QMN3)}. For (temporarily fixed) $\eps_2$ we choose

\begin{itemize}
	\item $n$ such that $\chi_0^n(\Ffamily) \leq \chi_0(\Ffamily) + \frac{\eps_2}{2}, \ \chi_0^n(\Gfamily) \leq \chi_0(\Gfamily) + \frac{\eps_2}{2}$ and 
			
	$$\min_{f \in \Ffamily} \ \inf_{g \in \Gfamily} \ d_{C(S_n,E)}(f,g) > 0$$
			
	\noindent
	In order to do that, pick $f_1 \in \Ffamily$ and choose $n \in \naturals$ such that 
			
	$$\inf_{g \in \Gfamily} \ d_{C(S_n,E)}(f_1,g) > 0$$
			
	\noindent
	This can be done, because otherwise $f_1 \in \Ffamily \cap \overline{\Gfamily}$, which is a contradiction. Then we increase (if necessary) $n$ so that 
			
			$$\inf_{g \in \Gfamily} \ d_{C(S_n,E)}(f_2,g) > 0$$
			
	\noindent
	where $f_2 \in \Ffamily\backslash\{f_1\}$. This step is repeated (finitely many times) until we exhaust all functions in $\Ffamily$.
	\item $\eps_1 > 0$ such that $\chi^n(\Ffamily,\eps_1) \leq \chi_0^n(\Ffamily) + \frac{\eps_2}{2}$ and $\chi^n(\Gfamily,\eps_1) \leq \chi_0^n(\Gfamily) + \frac{\eps_2}{2}$
\end{itemize}
 
\noindent
Define the function $\Psi : \Ffamily \rightarrow \reals_+$ by 
		
$$\Psi(f) := \inf_{g \in \Gfamily} \ d_{C(S_n,E)}(f,g)$$
		
\noindent
Intuitively, it measures the distance between the function $f$ and the set $\Gfamily$. By finiteness of $\Ffamily$, we obtain that there exists $f_{\ast} \in \Ffamily$ such that 
		
$$\min_{f \in \Ffamily} \ \Psi(f) = \Psi(f_{\ast})$$
		
\noindent
which is positive by our choice of $n$. If needed, we decrease $\eps_1$ so that $\eps_1 \leq \frac{1}{2}\Psi(f_{\ast})$. We have that 

\begin{gather}
\forall_{f,g \in C^b(X,E)} \ d_{C(S_n,E)}(f,g) \leq \eps_1 \ \Longrightarrow \ \bigg(f,g \in \Ffamily \ \vee \ f,g \in \Gfamily\bigg)
\label{splitab}
\end{gather}
		
\noindent
This implies that 
		
\begin{gather*}
\chi_0(\Ffamily \cup \Gfamily) \leq \sup\bigg\{d_{C^b(X,E)}(f,g) \ : \ f,g \in \Ffamily \cup \Gfamily, \ d_{C(S_n,E)}(f,g) \leq \eps_1 \bigg\} \\
\stackrel{(\ref{splitab})}{=} \max{\bigg(\chi^n(\Ffamily,\eps_1), \chi^n(\Gfamily,\eps_1)\bigg)} \leq \max{\bigg(\chi_0^n(\Ffamily), \chi_0^n(\Gfamily)\bigg)} + \frac{\eps_2}{2} \leq \max{\bigg(\chi_0(\Ffamily),\chi_0(\Gfamily)\bigg)} + \eps_2
\end{gather*}
	
\noindent
By theorem \ref{AAforXE}, since $\Ffamily$ is compact, we have $\chi_0(\Ffamily) = 0$. We end the proof by noting that $\eps_2$ can be arbitrarily small. 
\end{pro}

Finally, we are able to introduce the quasimeasure of noncompactness on $C^b(X,E)$. Piecing together all the considerations up to this point, we obtain the following theorem.

\begin{thm}(quasimeasure of noncompactness on $C^b(X,E)$)\\
For every $\Ffamily \in \bounded(C^b(X,E))$ we define

\begin{gather*}
\Omega(\Ffamily) := \eta(\Ffamily) + \omega_0(\Ffamily) + \chi_0(\Ffamily)
\end{gather*}

\noindent
The function $\Omega : \bounded(C^b(X,E)) \rightarrow \reals_+$ is a quasimeasure of noncompactness on $C^b(X,E)$. Moreover, for every $\Ffamily \in \bounded(C^b(X,E))$ we have 

\begin{gather}
\Omega(\Ffamily) = 0 \ \Longrightarrow \ \Omega(\conv(\Ffamily)) = 0
\label{Mproperty}
\end{gather}
\label{quasioncbxe}
\end{thm}
\begin{pro}

By lemmas \ref{pointwisecompactmeasurelemma}, \ref{equicontinuitymeasure}, \ref{chilemma} and theorem \ref{AAforXE} we obtain properties \textbf{(QMN1)} - \textbf{(QMN5)}. The property (\ref{Mproperty}) is a consequence of Mazur theorem, which can be found in \cite{Aliprantis} (theorem 5.35 on page 185). 
\end{pro}

\section{Application of quasimeasure of noncompactness in integral operators}
\label{applications}

In order to apply the machinery of quasimeasure of noncompactness in integral operators, we will need a Darbo-type theorem. However, the classical version (theorem 5.30, page 178 in \cite{BanasMursaleen}) will not be sufficient due to the fact that quasimeasure of noncompactness need not be invariant under convex hull (although it satisfies (\ref{Mproperty})), which is a fundamental property of measures of noncompactness. We are able to resolve this inconvenience by controlling the convexity with the aid of the measure of nonconvexity, which we recalled as definition \ref{mofnonconvexity}. 

The theorem we present below is stated in terms of $Y$ being any Banach space and $\Omega$ being any quasimeasure of noncompactness. However, in what follows, we focus on $Y = C^b(X,E)$ and the quasimeasure of noncompactness constructed in theorem \ref{quasioncbxe}.

\begin{thm}
Let $\Omega : \bounded(Y) \rightarrow \reals_+$ be a quasimeasure of noncompactness, $\kappa : \bounded(Y) \rightarrow \reals_+$ be a measure of nonconvexity and $C \subset Y$ be a nonempty, bounded and closed subset. If $\Phi : C \rightarrow C$ is a continuous function such that  

\begin{description}
	\item[\hspace{0.4cm} (D)] for all $A \subset C$ we have 
	
	$$\Omega(\Phi(A)) \leq \varphi_D(\Omega(A))$$
	
	\noindent
	where $\varphi_D : \reals_+ \rightarrow \reals_+$ is a nondecreasing function such that $\lim_{n\rightarrow \infty} \varphi_D^{(n)}(t) = 0$ for every $t \geq 0$.
	\item[\hspace{0.4cm} (E)] for all $A \subset C$ we have 
	
	$$\kappa(\Phi(A)) \leq \varphi_E(\kappa(A))$$
	
	\noindent
	where $\varphi_E : \reals_+ \rightarrow \reals_+$ is a nondecreasing function such that $\lim_{n\rightarrow \infty} \varphi_E^{(n)}(t) = 0$ for every $t \geq 0$.
\end{description}

\noindent
Then $\Phi$ has at least one fixed point in $C$. 
\label{Darbotheorem} 
\end{thm}
\begin{pro}

Let $C_1 := C$ and define $C_{n+1} := \overline{\Phi(C_n)}$. Closedness of $C_n$ for every $n\in\naturals$ is obvious and by induction, $(C_n)_{n\in\naturals}$ is a descending sequence. Defining $D := \bigcap_{n\in\naturals} C_n$ we see that it is a closed set (possibly empty at this stage of the proof).

By \textbf{(D)}, we have that

\begin{gather}
\Omega(C_{n+1}) = \Omega(\Phi(C_n)) \leq \varphi_D(\Omega(C_n)) \leq \varphi_D^{(2)}(\Omega(C_{n-1})) \leq \ldots \leq \varphi_D^{(n)}(\Omega(C))
\label{maindarboinequality}
\end{gather}

\noindent
which implies $\Omega(C_n) \rightarrow 0$ for $n\rightarrow \infty$. Hence $\Omega(D) = 0$, which means that $D$ is compact (in view of proved closedness). However, it may still be an empty set, which in the next step we prove is not the case. 

Choose a sequence $(x_n)_{n\in\naturals}\subset Y$ such that $x_n\in C_n$ and $x_n \neq x_m$ for every $n,m\in\naturals$. If it turns out that choosing such a sequence is impossible, then $D$ is a singleton. For a fixed $k \in \naturals$ we have by \textbf{(QMN5)} that 

$$\Omega((x_n)_{n\in\naturals}) = \Omega\bigg((x_n)_{n=1}^k \cup (x_{k+n})_{n\in\naturals}\bigg) \leq \Omega(C_{k+1})$$

\noindent
As $k\rightarrow\infty$, we conclude that $\Omega((x_n)_{n\in\naturals}) = 0$. By \textbf{(QMN1)} we can choose a convergent subsequence of $(x_n)_{n\in\naturals}$. The limit of this subsequence is in $D$, proving that $D \neq \emptyset$.

Similarly to (\ref{maindarboinequality}), by using \textbf{(E)}, we obtain $\kappa(C_n) \rightarrow 0$ for $n\rightarrow \infty$. By theorem 1 in \cite{Lesniak} we know that $d_H(C_n,D) \rightarrow 0$ and since 

$$|\kappa(C_n) - \kappa(D)| \leq 2d_H(C_n,D)$$

\noindent
we conclude that $\kappa(D) = 0$. Finally, we established that $D$ is nonempty, compact, convex and $\Phi(D)\subset D$. By Brouwer theorem, $\Phi$ has a fixed point. 
\end{pro}

In the sequel, we are preoccupied with obtaining a fixed point for Hammerstein operator. For further considerations, we need to introduce the concept of Carath\'{e}odory kernel. Our definition is inspired by definition 4.50 in \cite{Aliprantis} on page 153.

\begin{defin}(Carath\'{e}odory kernel)\\
We say that the function $K : X \times X \rightarrow \reals_+$ is the Carath\'{e}odory kernel if

\begin{description}
	\item[\hspace{0.4cm} (Car1)] $K(x,\cdot)$ is measurable for every $x \in X$
	\item[\hspace{0.4cm} (Car2)] $K(\cdot,y)$ is continuous for a.e. $y \in X$
	\item[\hspace{0.4cm} (Car3)] for every $x_{\ast} \in X$ there exists an open neighbourhood $U_x$ and a function $D_{x_{\ast}} \in L^1(X)$ such that for every $x \in U_x$ and a.e. $y \in X$ we have
	
	$$|K(x,y)| \leq D_x(y)$$
	
	\item[\hspace{0.4cm} (Car4)] kernel $K$ satisfies
	
	$$\sup_{x \in X} \ \int_X \ |K(x,y)| \ d\lambda(y) < \infty$$
\end{description}

\end{defin}

In the next theorem, we work on the cone

$$\cone(r,c) = \bigg\{f\in C^b(X,\reals_+) \ : \ \inf_{\|x\|\leq p} \ f(x) \geq c\|f\| \bigg\}$$

\noindent
where $r > 0$ and $c \in (0,1)$.

\begin{thm}
Let $K : X \times X \rightarrow \reals_+$ be a Carath\'{e}odory kernel and $N : X \times \reals_+ \rightarrow \reals_+$ be a nonlinearity of Hammerstein operator $\Hammer : C^b(X,\reals_+) \rightarrow \cone(r,c)$ defined by

$$(\Hammer f)(x) := \int_X \ K(x,y)N(y,f(y)) \ d\lambda(y)$$

\noindent
We assume that :

\begin{description}
	\item[\hspace{0.4cm} (N1)] for every $z_{\ast} \in \reals_+$ the function $N(\cdot,z)$ is measurable
	\item[\hspace{0.4cm} (N2)] for every $z_{\ast} \in \reals_+$ and $\eps > 0$ there exists $\delta > 0$ such that for a.e. $y \in X$ and $z \in \reals_+$ we have
	
		$$\|z - z_{\ast}\|\leq \delta \ \Longrightarrow \ |N(y,z)-N(y,z_{\ast})| \leq \eps$$
	
	\item[\hspace{0.4cm} (K1)] for a.e. $y \in X$ we have $K(\cdot, y) \in \cone(r,c)$
	
	\item[\hspace{0.4cm} (K2)] there exists a continuous function $\zeta : X \times \reals_+ \rightarrow \reals_+$ nondecreasing in the second variable, such that for every $y \in X, \ z\in \reals_+$ we have $N(y,z) \leq \zeta(y,z)$ and there is a $R \in (0,\infty)$ satisfying 
	
	$$\sup_{x \in X} \ \int_X \ K(x,y) \ \zeta(y,R) \ d\lambda(y) = R$$
	
	\item[\hspace{0.4cm} (H)] there exists $q \in (0,1)$ such that for every $\Ffamily \subset B(0,R) \cap \cone(r,c)$ we have 
	
	$$\chi_0(\Hammer(\Ffamily)) \leq q \ \chi_0(\Ffamily)$$
\end{description}

\noindent
Then $\Hammer$ has a fixed point.
\label{fixedpointdarbo} 
\end{thm}
\begin{pro}

In \cite{KrukowskiPrzeradzki} we find that the Carath\'{e}odory kernel together with assumptions \textbf{(N1)}-\textbf{(N2)} are sufficient for the Hammerstein operator to be continuous and map $C^b(X,\reals_+)$ into itself. Moreover, the calculation 

\begin{equation}
	\begin{split}
		&\inf_{\|x\| \leq r} \ \int_X \ K(x,y)N(y,f(y)) \ d\lambda(y) \geq \int_X \ \inf_{\|x\| \leq r} \ K(x,y)N(y,f(y)) \ d\lambda(y) \\
		&\stackrel{\textbf{(K1)}}{\geq} \int_X \ c \sup_{x \in X} \ K(x,y)N(y,f(y)) \ d\lambda(y) \geq c \sup_{x \in X} \ \int_X \ K(x,y)N(y,f(y)) \ d\lambda(y)
	\end{split}
	\label{whyhammersteinisincone}
\end{equation}

\noindent
proves that the image of Hammerstein operator is in fact in $\cone(r,c)$.

Denote $C = B(0,R) \cap \cone(r,c)$ and observe that for every $f \in C$ we have 

\begin{equation}
	\begin{split}
		\|\Hammer &f\|_{C^b(X,\reals_+)} \stackrel{\textbf{(K2)}}{\leq} \sup_{x \in X} \ \int_X \ K(x,y) \ \zeta(y,|f(y)|) \ d\lambda(y) \\
		&\leq \sup_{x \in X} \ \int_X \ K(x,y) \ \zeta(y,R) d\lambda(y) \stackrel{\textbf{(K2)}}{=} R 
	\end{split}
	\label{inequalitywithphi}
\end{equation}

\noindent
By (\ref{whyhammersteinisincone}) and (\ref{inequalitywithphi}), we conclude that $\Hammer(C) \subset C$.

In order to verify \textbf{(D)} in theorem \ref{Darbotheorem}, we observe that for all $\Ffamily \subset C$ we have $\eta(\Ffamily) = 0$. Moreover, due to the fact that we are working with Carath\'{e}odory kernel we have $\omega_0(\Hammer(\Ffamily)) = 0$ (proof in \cite{KrukowskiPrzeradzki}). By \textbf{(H)} we obtain \textbf{(D)}.

Lastly, observe that we do not need to verify condition \textbf{(E)}. This is because $C$ is convex as the intersection of two convex sets. The proof of theorem \ref{Darbotheorem} works now without invoking the measure of nonconvexity. Finally, we conclude the existence of a fixed point.
\end{pro}

As a final note, let us explain why theorem \ref{fixedpointdarbo} may cause some problems in practice. The difficulty lies in assumption \textbf{(H)}. This condition is obviously satisfied, if the nonlinearity $N$ is a contraction in the second variable (the contracting constant $q$ is chosen uniformly with respect to the first variable). This is compatible with the spirit of Banach contraction principle. 

Another situation, when \textbf{(H)} is satisfied is when (additionally to assumptions in theorem \ref{fixedpointdarbo}) for every $\eps > 0$ and $x \in X, \ \|x\| = 1$ there exist $T_x > 0$ and $L_x \in L^1(X,\reals_+)$ such that for every $t \geq T_x$ we have
	
$$\int_X \ |K(tx,y)-L_x(y)| \ d\lambda(y) \leq \eps$$
	
\noindent
and moreover $\sup_{\|x\|=1} T_x < \infty$. In \cite{KrukowskiPrzeradzki}, Hammerstein operators with such kernels have been proved to be compact. Hence \textbf{(H)} is obviously true.

In view of these remarks, the author would like to point out the similarity between classical Darbo theorem and theorem \ref{fixedpointdarbo}. Darbo theorem is not a new technique for generating fixed point. At the heart of the proof lies a very powerful Schauder theorem, which is the true source of the fixed point. Darbo's result binds both Banach- and Schauder-like approaches, which is a beautiful example of unification in mathematics. The author believes the same holds true for theorem \ref{fixedpointdarbo}.

\section*{Acknowledgements}

I am particularly grateful for the assistance given by my supervisor Bogdan Przeradzki. His constant support was invaluable during my work on this paper.

I also wish to acknowledge the help provided by Dariusz Bugajewski, who suggested to improve the results by introducing the measure of nonconvexity. I highly appreciate his invitation to the seminar and fruitful discussion of my work.

\Addresses
\end{document}